\documentclass[12pt]{article}

\usepackage{amsmath,amsthm,indentfirst}
\usepackage{amssymb}
\usepackage{amsfonts}
\usepackage{amscd}
\usepackage[latin1]{inputenc}

\theoremstyle{plain}
\newtheorem*{maintheorem*}{Main Theorem}
\newtheorem*{thm*}{Theorem}
\newtheorem*{thma*}{Theorem A}
\newtheorem*{thmaa*}{Theorem A'}
\newtheorem*{thmb*}{Theorem B}
\newtheorem*{thmo*}{Theorem 1}
\newtheorem*{thmc*}{Theorem C}
\newtheorem*{thmd*}{Theorem D}
\newtheorem*{thmf*}{Theorem 4.1}
\newtheorem*{conjecture*}{Conjecture}
\newtheorem*{proposition*}{Proposition}
\newtheorem{thm}{Theorem}

\newtheorem{lem}[thm]{Lemma}
\newtheorem{proposition}[thm]{Proposition}
\newtheorem{example}[thm]{Example}

\theoremstyle{definition}

\newtheorem*{proofc*}{Proof of Theorem C}

\newtheorem{conjecture}[thm]{Conjecture}

\newtheorem{definition}[thm]{Definition}

\newtheorem{remark}[thm]{Remark}

\title{Non-commutative Henselian Rings} 
\author{ Masood Aryapoor}
 
\begin{document}
 \maketitle
\begin{abstract}
\noindent
Non-commutative Henselian rings are defined and it is shown that a local ring  which is complete and separated in the topology defined by its maximal ideal is Henselian provided that it is almost commutative.
\end{abstract}
 
 We define non-commutative Henselian rings and give some examples of them.
 Here, all rings are assumed to be unitary.
Let us start with a definition,\\
\begin{definition}
A (possibly non-commutative) ring $A$ is called local if all the non-invertible
elements form an (two-sided) ideal which we denote by $m$.
\end{definition}
If $A$ is a local ring, then  $k=A/m$ is a skew field,
called the residue field. We denote the reduction map $A\to k$ by ($a\to \bar{a}$).
For a brief introduction to local rings consult ~\cite{Lam}, Chapter 7.\\
Let $A[x]$ be the ring of polynomials over $A$ where the
indeterminate $x$ commutes with elements of $A$. Commutative Henselian rings are defined as follows,
\begin{definition}
Let $A$ be a commutative local ring with the maximal ideal $m$ and residue field $k$. $A$ is called Henselian if  for
every   polynomial $f(x)=x^n+a_{n-1}x^{n-1}+\cdots+a_1x+a_0\in A[x]$ such that
$\overline{f(x)}=f_1(x)f_2(x)$ for some relatively prime
monic polynomials  $f_i(x)\in k[x]$ then there are  unique monic polynomials
$F_i(x)\in A[x]$ such that $f(x)=F_1(x)F_2(x)$ and $\overline{F_i(x)}=f_i(x)$.
\end{definition}
See ~\cite{Ray} for a detailed discussion of commutative Henselian rings.\\
The above definition makes sense as long as $k$, the residue field, is commutative. 
Therefore we have the following definition,
\begin{definition}
Let $A$ be a (possibly non-commutative) local ring  with the maximal ideal $m$ and residue field $k$.
Moreover assume that $k$ is commutative. Then $A$ is called Henselian if  for
every   polynomial $f(x)=x^n+a_{n-1}x^{n-1}+\cdots+a_1x+a_0\in A[x]$ such that
$\overline{f(x)}=f_1(x)f_2(x)$ for some relatively prime
monic polynomials  $f_i(x)\in k[x]$ then there are  unique monic polynomials
$F_i(x)\in A[x]$ such that $f(x)=F_1(x)F_2(x)$ and $\overline{F_i(x)}=f_i(x)$.
 \end{definition}
 It is well-known that every commutative local ring $A$ which is complete and  separated in the $m$-adic topology 
 is Henselian. This is not true for non-commutative local rings which are complete and separated  in the topology defined by the maximal ideal. However, it holds if the local ring has an extra property which we explain in what follows.\\
 To each local ring one can associate an associative ring as
follows,
\begin{definition}
Let $A$ be a local ring with the maximal ideal $m$. Then $gr(A)=\frac{A}{m}\bigoplus \frac{m}{m^2}\bigoplus \cdots$ is defined to be the graded associated ring coming from
the filtration  $\cdots\subset m^{n+1}\subset m^{n} \subset \cdots \subset m
\subset A$. $A$ is called almost commutative if $gr(A)$ is
commutative.
\end{definition}
For basic facts regarding $gr(A)$ see [Lang].\\
Clearly if $A$ is almost commutative, then $k$ is commutative.
The main theorem is,
\begin{thm}
Let $A$ be an almost commutative local ring such that $A$ is both
separated, i.e. $\bigcap m^{n}=\{0\}$,  and complete in the $m$-adic topology. Then $A$ is a
Henselian ring.
\end{thm}
\begin{proof}
Basically, the same proof of Hensel's lemma works.
Let $f(x)=x^n+a_{n-1}x^{n-1}+\cdots+a_1x+a_0\in A[x]$  such that $\overline{f(x)}=f_1(x)f_2(x)$ for some relatively prime monic
polynomials  $f_i(x)\in k[x]$.   We will inductively construct a   sequence of monic polynomials $\{ F_{1,r}(x) \}$ and $\{ F_{2,r}(x) \}$ in $A[x]$ such that $\overline{ F_{1,r}(x)}=f_1(x)$,  $\overline{ F_{2,r}(x)}=f_2(x)$, $F_{1,r+1}(x)-F_{1,r}(x) \in m^{r}[x]$, $F_{2,r+1}(x)-F_{2,r}(x) \in m^{r}[x]$ and $f(x)-F_{1,r}(x)F_{2,r}(x)\in m^{r}[x]$. Clearly this proves the existence part.\\
It is easy to find $F_{1,1}(x)$ and $F_{2,1}(x)$. Having defined $F_{1,r}(x)$ and $F_{2,r}(x)$, we define $F_{1,r+1}(x)$ and $F_{2,r+1}(x)$ as follows. Writing $F_{1,r+1}(x)=F_{1,r}(x)+G_1(x)$ and $F_{2,r+1}(x)=F_{2,r}(x)+G_2(x)$, finding $F_{1,r+1}$ and $F_{2,r+1}$ is equivalent to finding $G_1(x)$ and $G_2(x)$ in $m^{r}[x]$ such that 
$deg(G_1(x))<deg(f_1(x))$, $deg(G_2(x))<deg(f_2(x))$ and $$f(x)-F_{1,r}(x)F_{2,r}(x)-G_1(x)F_{2,r}(x)-F_{1,r}(x)G_2(x)\in m^{r+1}[x].$$ By abuse of notations this is the same as finding  $G_1(x)$ and $G_2(x)$ in $m^{r}[x]$ such that 
$deg(G_1(x))<deg(f_1(x))$, $deg(G_2(x))<deg(f_2(x))$ and $f(x)-F_{1,r}(x)F_{2,r}(x)-G_1(x)F_{2,r}(x)-F_{1,r}(x)G_2(x)=0$ in  $m^{r}/m^{r+1}[x]$. Considering $m^{r}/m^{r+1}$ as a vector space over $k=A/m$ and using the fact that $A$ is almost commutative, one can see that this is the same as finding $G_1(x)$ and $G_2(x)$ in $m^{r}[x]$ such that 
$deg(G_1(x))<deg(f_1(x))$, $deg(G_2(x))<deg(f_2(x) )$ and $(f(x)-F_{1,r}(x)F_{2,r}(x))-f_2(x)G_1(x)-f_1(x)G_2(x)=0$ in  $m^{r}/m^{r+1}[x]$. This is possible because $f_1(x)$ and $f_2(x)$ are relatively prime.\\
The uniqueness part follows from the facts that  $f_1(x)$ and $f_2(x)$ are relatively prime and $A$ is separated in the $m$-adic topology.
 \end{proof}
In the commutative case, one can use Hensel's lemma to find roots of polynomials. Next we show this connection in the non-commutative case.\\
Let $A[x]$ be the ring of polynomials over $A$ where the
indeterminate $x$ commutes with elements of $A$. So every element
of $f(x)\in A[x]$ can be written uniquely as
$f(x)=a_nx^n+\cdots+a_1x+a_0$ with $a_i\in A$. One can consider
$f(x)$ as a function on $A$ as follows,
$f(a):=a_na^n+\cdots+a_1a+a_0$ for $a\in A$. 
\begin{definition}
An element $a\in A$ is
called a (right) root of $f(x)=a_nx^n+\cdots+a_1x+a_0$ if $f(a)=0$.
\end{definition}
 We have the following proposition,
 \begin{proposition}
 An element $a\in A$ is a root of $f(x)=a_nx^n+\cdots+a_1x+a_0\in A[x]$ if and only if $f(x)=g(x)(x-a)$ for some $g(x)\in A[x]$.
 \end{proposition}
  To see the proof and basic facts regarding right and left roots, see ~\cite{Lam}, Chapter 5.\\
  Theorem 5 together with the above proposition imply that,
  \begin{thm}
  Let $A$ be a Henselian ring. Suppose that $f(x)=x^n+a_{n-1}x^{n-1}+\cdots+a_1x+a_0\in A[x]$ is a monic polynomial such that $\overline{f(x)}$ has a simple root $r\in k$. Then $f(x)$ has a unique root $a\in A$ such that $\bar{a}=r$.
 \end{thm}
In the commutative case, a local ring $A$ is Henselian if and only if every finite $A$-algebra is isomorphic to a product of local rings (See [Ray]). In the non-commutative case we can give a similar criterion for Henselian rings in terms of some conditions on some modules over $A$.\\
We begin with a few definitions.
\begin{definition}
Let $A$ be a ring and $M$ a (left) $A$-module. We say that $M$ is local if it has a unique maximal submodule. $M$ is called semi-local if $M=M_1\bigoplus \cdots \bigoplus M_k$ where $M_i$'s are local. It is called indecomposable if it cannot be written as $M=M_{1}\bigoplus M_{2}$, where $M_{i}$'s are nonzero submodule of $M$. It is called strongly indecomposable if $End_A(M)$ is a local ring.
\end{definition} 
One has the following theorem.\\
\begin{thm}{(Krull-Schmidt-Azumaya)}

Suppose that the $A$-module $M$ has the following decompositions into submodules,
$$M=M_1\bigoplus \cdots \bigoplus M_r\simeq N_1\bigoplus \cdots \bigoplus N_s,$$
where $M_i$'s are indecomposable and $N_i$'s are strongly indecomposable. Then $r=s$ and after a reindexing we have $M_i\simeq N_i$.
\end{thm}
For a proof see [Lam], chapter seven.\\
From now on, we suppose that $A$ is a local ring as before. Let $M$ be an $A$-module.   Set $\bar{M}=\frac{M}{mM}$ which is a $k$-module. 
We need a few lemmas. 
\begin{lem}
Let $M$ be an $A[x]$-module which is a finitely generated $A$-module. Then, $M$ is a local $A[x]$-module if and only if $\frac{M}{mM}$ is a local $k[x]$-module.
\end{lem}
\begin{proof}
By Nakayama's lemma, every maximal submodule of $M$ contains $mM$.
\end{proof}
\begin{lem}
Let $M,N$ be finitely generated $A$-modules. Let $\alpha:M\to N$ be an $A$-module homomorphism and $\overline{\alpha}:\bar{M}\to \bar{N}$ be the induced $k$-linear map. If $ker(\overline{\alpha})\neq 0$ and $\bar{\alpha}$ is onto, then $ker(\alpha)\neq 0$.
\end{lem}
\begin{proof}
Suppose that $v_1,...,v_n$ are elements of $M$ such that $\overline{\alpha}(\overline{v_1}),...,\overline{\alpha}(\overline{v_n})$ form a basis for $\bar{N}$ over $k$. Then by Nakayama's lemma we have that $\alpha(v_1),...,\alpha(v_n)$ generate $N$ as an $A$-module. Since $ker(\bar{\alpha})\neq 0$, $\overline{v_1},...,\overline{v_n}$ do not generate $\bar{M}$. So $v_1,...,v_n$ do not generate $M$ which follows that $ker(\alpha)\neq 0$.
\end{proof}
We also need the following lemma,
\begin{lem}
Let $A$  be a local ring whose residue field $k$ is commutative.
Suppose that $p,q\in A[x]$ are polynomials of degrees $r,s$
respectively and $p$ is monic. If $A[x]p+A[x]q=A[x]$, then there
are polynomials $p_1,q_1\in A[x]$ such that $deg(p_1)=deg(q)$,
$deg(p)=deg(q_1)$, $p_1p=q_1q$ and $q_1$ is monic.
 
\end{lem}
\begin{proof}
Let $\alpha:A^{s+1}\bigoplus A^{r+1}\to A^{r+s+1}$ be the following map,
$$\alpha(a_0,a_1,...,a_s,b_0,b_1,...,b_r)=(\sum_{i=0}^{s}{a_ix^i})p-(\sum_{i=0}^{r}{b_ix^i})q.$$
Using lemma 12, we have that $ker(\alpha)\neq 0$. This shows that there are nonzero polynomials $p_1,q_1\in A[x]$ such that $deg(p_1)\leq deg(q)$, $deg(q_1)\leq deg(p)$, $p_1p=q_1q$. Since $\bar{p}$ and $\bar{q}$ are prime in $k[x]$ and $p$ is monic, we must have $deg(\bar{q_1})=deg(p)$, hence $deg(q_1)=deg(p)$ and $deg(p_1)=deg(q)$. Finally, it is clear that $q_1$ can be chosen to be monic. 

\end{proof}
\begin{remark}
If $p,q\in A[x]$ are polynomials such that $p$ is monic and  $A[x]p+A[x]q+m[x]=A[x]$ then $A[x]p+A[x]q=A[x]$. In fact we have that $M=\frac{A[x]}{A[x]p+A[x]q}$ is a finitely generated $A$-module and $mM=M$. So, by Nakayama's lemma, $M=0$.
\end{remark}
We have the following theorem,
\begin{thm}
Suppose that $A$ is a local ring whose residue field $k$ is commutative. Then the following are equivalent,\\
(1) $A$ is Henselian.\\
(2) For any monic polynomial $p\in A[x]$ the $A[x]$-module $M=\frac{A[x]}{A[x]p}$ is semi-local.\\
 \end{thm}

\begin{proof}
First we show that (1) implies (2). If $\bar{p}$ is a power of an irreducible polynomial in $k[x]$ then 
$\bar{M}=\frac{M}{mM}=\frac{k[x]}{(\bar{p})}$ is a local $k[x]$-module and by   lemma 11, $M$ is local.
Suppose $\bar{p}=f_1f_2$ where $f_1$ and $f_2$ are relatively prime polynomials of $k[x]$. By (1) we have $p=p_1p_2=q_2q_1$ where $p_i,q_i$ are monic polynomials in $A[x]$ such that $\overline{p_i}=\overline{q_i}=f_i$. This implies that $M\simeq \frac{A[x]}{A[x]p_2}\bigoplus \frac{A[x]}{A[x]q_1}$ because $A[x]p_2+A[x]q_1=A[x]$(above remark) and it is easy to see that $A[x]p_2\cap A[x]q_1=A[x]p$. Now we can use induction on $deg(p)$.\\
Conversely, let $p\in A[x]$ be a monic polynomial. Then we have $M=\frac{A[x]}{A[x]p}=M_1\bigoplus \cdots \bigoplus M_r$ where $M_i$'s are local. So we have $\bar{M}=\bar{M_1}\bigoplus \cdots \bigoplus \bar{M_r}$. On the other hand, if $\bar{p}=f_1\cdots f_s$ where $f_i$'s are powers of irreducible monic  polynomials in $k[x]$, then $\bar{M}\simeq \frac{k[x]}{(f_1)}\bigoplus \cdots \bigoplus  \frac{k[x]}{(f_s)}$. It is easy to see that  $ \frac{k[x]}{(f_i)}$'s  are strongly indecomposable as $k[x]$-modules and $\bar{M_i}$'s are local, in particular indecomposable. So by Krull-Schmidt-Azumaya theorem  $r=s$ and $\bar{M_{i}}\simeq \frac{k[x]}{(f_i)}$ possibly after a reindexing $M_{i}$'s. If $v_i\in M_i$ is the image of  $1\in A[x]$ then $(Av_i+Axv_i+\cdots+Ax^{n_{i}-1}v_i)+mM_i=M_i$  where $n_i$ is the degree of $f_i$.
By Nakayama's lemma  $(Av_i+Axv_i+\cdots+Ax^{n_{i}-1}v_i)=M_i$. Also $p_iv_i=0$ for some monic  polynomial $p_i$ of degree $n_i$ such that $\bar{p_i}=f_i$. By lemma 13, there is a monic polynomial  $p'=q_1q_2\cdots  q_r$ where $q_i$'s are monic polynomials and $\bar{q_i}=f_i$ and $p'\in A[x]p_i$ for each $i$. This implies that $p'\in A[x]p$. Since $deg(p)=deg(p')$ and they are monic we have $p'=p$.

\end{proof}

Finally we give some examples.
\begin{example}
Let $k$ be a  field with a derivation. The ring of Volterra operators $k[[\partial^{-1}]]$ is defined as follows(See [Lebedev] for more on Volterra operators). 
It is the set of formal series  $a_0+a_{1}\partial^{-1}+\cdots$ with $a_i \in k$ where $\partial^{n} a=\sum_{i=0}^{\infty}{{n\choose{i}} a^{(i)}\partial^{n-i}}$ for $n<0$. One can see that $k[[\partial^{-1}]]$ is a local ring with the maximal ideal $m=k[[\partial^{-1}]]\partial^{-1}$ which is both separated and complete in the
$m$-adic topology. Moreover $gr(k[[\partial^{-1}]])$ is isomorphic to $k[x]$ the ring of polynomials over $k$, hence commutative. So  $k[[\partial^{-1}]]$ is a Henselian ring.
 \end{example}

\begin{example}
 If $A$ is not almost  commutative but complete and separated in the $m$-adic  topology then there might not be any lifting of simple roots. Here is one example.
Let $k$ be a field and $\sigma$ an
automorphism of $k$. Let $A$ be the set of all series of the form
$a_0+a_1\tau+a_2\tau^2+\cdots$  where $a_i \in k$. One can make $A$
into a ring using the relation $\tau a=\sigma(a)\tau$ for $a\in
k$. Then $A$ is a local ring which is both separated and complete
in the $m$-adic topology and $A/m=k$ is commutative. However if
$\sigma$ is not the identity map then $gr(A)$ is isomorphic to the skew polynomial ring $k[x;\sigma]$, hence not commutative. 
Suppose $k=\mathbb{C}$ and $\sigma$ is the complex conjugation.  Consider the polynomial $f(x)=x^2+1+\tau$ in $A[x]$. Then
$\overline{f(x)}$ has  a simple root in $k$, namely  $\sqrt{-1}$. However $f(x)$ does not
have any root in $A$. Since if $a=a_0+a_1\tau+a_2\tau^2+\cdots$ is a root of $f(x)$ then we have 
$0=a^2+1+\tau=a_0^2+1+(a_0a_1+\overline{a_0}a_1+1)\tau+\cdots$. This implies that $a_0=\sqrt{-1}$ or $a_0=-\sqrt{-1}$. Therefor $a_0a_1+\overline{a_0}a_1+1=1\neq0$, a contradiction.

\end{example}
In the commutative case, for any local Noetherian ring $A$, there is a (unique) Henselian ring $A^{h}$, called  the Henselization of $A$, and a local homomorphism $i:A\to A^{h}$ with the following universal property,
given any local homomorphism $f$ from $A$ to some Henselian ring $B$ there is a unique local homomorphism $f^{h}:A^{h}\to B$ such that $f=f^{h}i$. \\
One can ask the same question in the non-commutative case. If $A$ is a local ring such that $gr(A)$ is commutative, then the completion of $A$ with respect to the $m$-adic topology is Henselian provided that it is separated. It is easy to see that the intersection of all local Henselian rings $H$ in the completion $\hat{A}$, with the maximal ideal $m_H$ such that $A\subset H\subset \hat{A}$ and $m_{\hat{A}}\cap H=m_H$, denoted by $\bar{A}$, is a Henselian local ring. In the commutative case it is not hard to see that $\bar{A}$ is the Henselization.
Therefore one might propose the following conjecture,
\begin{conjecture}
The Henselization exists for any almost commutative separated local ring $A$ and $A^h\simeq \bar{A}$. 
\end{conjecture}


 Department of Mathematics, Yale University, 10 Hillhouse Avenue, New Haven
CT 06520 USA, email: masood.aryapoor@yale.edu

 \end{document}